# COEFFICIENT ESTIMATES FOR INVERSES OF STARLIKE FUNCTIONS OF POSITIVE ORDER

**G. P. KAPOOR**
Department of Mathematics
Indian Institute of Technology Kanpur
Kanpur 208016, India
e-mail: gp@iitk.ac.in

and

**A. K. MISHRA***
Department of Mathematics
Berhampur University
Berhampur 760007, India
e-mail: akshayam2001@yahoo.co.in

**ABSTRACT:** In the present paper, the coefficient estimates are found for the class $S^{*-1}(\alpha)$ consisting of inverses of functions in the class of univalent starlike functions of order $\alpha$ in $D = \{z \in C : |z| < 1\}$. These estimates extend the work of Krzyz, Libera and Zlotkiewicz [12] who found sharp estimates on only first two coefficients for the functions in the class $S^{*-1}(\alpha)$. The coefficient estimates are also found for the class $\Sigma^{*-1}(\alpha)$, consisting of inverses of functions in the class $\Sigma^*(\alpha)$ of univalent starlike functions of order $\alpha$ in $V = \{z \in C : 1 < |z| < \infty\}$. The open problem of finding sharp coefficient estimates for functions in the class $\Sigma^*(\alpha)$ stands completely settled in the present work by our method developed here.

**AMS (MOS) Subject classification:** Primary *30C50*; Secondary *30C45*

*(\*) The present research of the author is partially supported through Grant No. 48/2/2003-R&D-II/1158, National Board for Higher Mathematics, Department of Atomic Energy, Government of India.*

## 1. INTRODUCTION

Let $A_0$ be the class of functions $f$, analytic in the unit disc $D = \{z \in C : |z| < 1\}$ and having the power series expansion,

$$f(z) = z + \sum_{n=2}^{\infty} a_n z^n. \tag{1.1}$$

The class of univalent functions in $A_0$ is denoted by $S$. A function $f \in A_0$ is said to be in the class $S^*(\alpha)$, $0 \leq \alpha < 1$, of starlike functions of order $\alpha$ if, for $z \in D$,

$$\operatorname{Re}\left(z\frac{f'(z)}{f(z)}\right) > \alpha, \tag{1.2}$$

It is known that $S^*(\alpha) \subseteq S^*(0) \equiv S^* \subset S$ for $0 \leq \alpha < 1$ [5, p.51]. The class of functions

$$g(z) = z + b_0 + \frac{b_1}{z} + \frac{b_2}{z^2} + \cdots \tag{1.3}$$

that are analytic and univalent in $V = \{z \in C : 1 < |z| < \infty\}$ is denoted by $\Sigma$ and the class of starlike functions of order $\alpha$, $0 \leq \alpha < 1$, in $V$ is denoted by $\Sigma^*(\alpha)$, i.e. a function $g \in \Sigma^*(\alpha)$ if and only if $g \in \Sigma$ satisfies $\operatorname{Re}\left(\frac{z\,g'(z)}{g(z)}\right) > \alpha$ for $z \in V$. We are primarily concerned here with the investigation of sharp coefficient estimates for the inverse functions in the above classes. Let $S^{-1}$ be the class of inverse functions $f^{-1}$ of functions $f \in S$ with the Taylor series expansion

$$f^{-1}(w) = w + \sum_{n=2}^{\infty} A_n w^n \tag{1.4}$$

in some disc $|w| < r_0(f)$ and $\Sigma^{*-1}$ be the class of inverse functions $g^{-1}$ of functions $g \in \Sigma^*$ with the series expansion

$$g^{-1}(w) = w + B_0 + \frac{B_1}{w} + \frac{B_2}{w^2} + \cdots \tag{1.5}$$

in some neighbourhood of infinity. The classes $S^{*-1}(\alpha)$ and $\Sigma^{*-1}(\alpha)$ are defined analogously.

The coefficient estimate problem for the class $S$, known as the Bieberbach conjecture [2], is settled by De-Branges [4], who proved that for a function $f(z) = z + \sum_{n=2}^{\infty} a_n z^n$ in the class $S$, $|a_n| \leq n$, for all $n = 2, 3, \ldots$, with equality only for the rotations of the Koebe function

$$K_0(z) = \frac{z}{(1-z)^2} \tag{1.6}$$

Loewner, using his parametric method [15; also see 7, p.222] proved that if $f^{-1}$, given by (1.4), is in the class $S^{-1}$ or $S^{*-1}$, the sharp estimate

$$|A_n| \leq \frac{\Gamma(2n+1)}{\Gamma(n+2)\,\Gamma(n+1)}, \qquad n = 2, 3, \ldots$$

holds, $K_0^{-1}$ being the extremal function for all $n$ in the above inequality. The above coefficient estimate problem for the classes $S^{-1}$ and $S^{*-1}$ is also investigated by the methods developed by Shaeffer and Spencer [17], FitzGerald [6], Baernstein [1], Poole [16] and others. Prior to de-Branges result on the sharp coefficient bounds for the whole class $S$, the coefficient estimate problem was established for several subclasses of $S$, e.g. the classes of convex functions, starlike functions of order $\alpha$, close-to-convex functions, normalized integrals of functions with positive real part etc. However, in contrast, although for the class of inverse functions the coefficient problem for the whole classes $S^{-1}$ and $S^{*-1}$ had been completely solved as early as in 1923 [15], only few results are known on the sharp coefficient estimates for the inverse of functions in the above subclasses [9, 11, 14, 18]. In certain cases, the coefficients of the inverse of the functions in some of these subclasses show unexpected behaviour. For example, it is known that if $f$ is a univalent convex function and $f^{-1}$ is given by (1.4), then $|A_n| \leq 1$ for $n = 1,...,8$ [3,13] and equality holds for the inverse of the function $K_{1/2}(z) = z/(1-z)$; while $|A_{10}| > 1$ [10] and the exact bounds on $|A_9|$ and $|A_n|$ for $n > 10$ are still unknown.

Krzyz, Libera and Zlotkiewicz [12] showed that if $f^{-1}$, given by (1.4), is in $S^{*-1}(\alpha)$, then
$$|A_2| \leq 2(1-\alpha) \tag{1.7}$$
and
$$|A_3| \leq \begin{cases} (1-\alpha)(5-6\alpha), & 0 \leq \alpha \leq \tfrac{2}{3} & (1.8) \\ (1-\alpha), & \tfrac{2}{3} \leq \alpha < 1 & (1.9) \end{cases}$$

The estimates (1.7) and (1.8) are sharp for the function $K_\alpha^{-1}$ and (1.9) is sharp for the function $K_{\alpha,2}^{-1}$, where
$$K_\alpha(z) = \frac{z}{(1-z)^{2(1-\alpha)}} \quad \text{and} \quad K_{\alpha,n}(z) = \sqrt[n]{K_\alpha(z^n)} \quad, n = 2,3,... \tag{1.10}$$

The determination of sharp estimates on $|A_n|$ for $n \geq 4$ is hitherto an open problem for the class $S^{*-1}(\alpha)$.

In the present paper, the estimates on $|A_n|$, $n \geq 2$, for functions in the class $S^{*-1}(\alpha)$ are found. These estimates are sharp for each $n$ when $\alpha \in [0, \tfrac{2}{n}) \cup [\tfrac{n-1}{n}, 1)$ extending the work of Krzyz, Libera and Zlotkiewicz [12], who found sharp estimates given by (1.7), (1.8), and (1.9) on only first two coefficients for the functions in the class $S^{*-1}(\alpha)$. More specifically, it is shown that for each $n \geq 3$, the functions $K_\alpha^{-1}$ and $K_{\alpha,n}^{-1}$, given by (1.10), are extremal for $0 \leq \alpha < \tfrac{2}{n}$ and $\tfrac{n-1}{n} \leq \alpha < 1$ respectively; at variance to the existence of a single extremal function $K_0^{-1}$ for the whole class $S^{*-1}$. We also prove that, for the functions $g^{-1} \in \Sigma^{*-1}(\alpha)$, given by (1.5), the sharp estimate $|B_n| \leq 2(1-\alpha)$, $\tfrac{n-1}{n} \leq \alpha < 1$, holds for all $n = 1,2,\cdots$. The open problem of finding sharp coefficient estimates for functions in the class $\Sigma^*(\alpha)$ stands completely settled in the present work by our method developed here.

## 2. SOME AUXILIARY RESULTS

Let $\Omega$ be the class of functions

$$f(z) = \sum_{v=1}^{\infty} a_v z^v, \qquad a_1 \neq 0 \tag{2.1}$$

analytic in $|z| < \rho$ for some $\rho > 0$. For a fixed $n \in N$ and $f \in \Omega$, let

$$\frac{1}{(f(z))^n} = \frac{a_{-n}^{(-n)}}{z^n} + \frac{a_{-n+1}^{(-n)}}{z^{n-1}} + \cdots + \frac{a_{-1}^{(-n)}}{z} + \sum_{p=0}^{\infty} a_p^{(-n)} z^p = \sum_{\vartheta=0}^{\infty} a_{-n+\vartheta}^{(-n)} z^{-n+\vartheta}, \quad |z| < \rho \tag{2.2}$$

Note that $a_{-n}^{(-n)} = 1$, if $a_1 = 1$.

We need the following lemmas in the sequel for our main results:

**Lemma1.** *Let $0 \leq \alpha < 1$ and $n \in N$ be fixed. Then,*

$$4n(1-\alpha)\left[n(1-\alpha) + \sum_{m=1}^{\vartheta-1}\{n(1-\alpha)-m\}\left(\prod_{j=0}^{m-1}\frac{2n(1-\alpha)-j}{j+1}\right)^2\right]$$

$$= \frac{1}{((\vartheta-1)!)^2} \prod_{j=0}^{\vartheta-1}(2n(1-\alpha)-j)^2 \tag{2.3}$$

**Proof.** For $\vartheta = 1$, (2.3) trivially holds. Assume that

$$4n(1-\alpha)\left[n(1-\alpha) + \sum_{m=1}^{\vartheta-2}\{n(1-\alpha)-m\}\left(\prod_{j=0}^{m-1}\frac{2n(1-\alpha)-j}{j+1}\right)^2\right] = \frac{1}{((\vartheta-2)!)^2}\prod_{j=0}^{\vartheta-2}(2n(1-\alpha)-j)^2.$$

Then,

$$4n(1-\alpha)\left[n(1-\alpha) + \sum_{m=1}^{\vartheta-1}\{n(1-\alpha)-m\}\left(\prod_{j=0}^{m-1}\frac{2n(1-\alpha)-j}{j+1}\right)^2\right]$$

$$= \frac{1}{((\vartheta-2)!)^2}\prod_{j=0}^{\vartheta-2}(2n(1-\alpha)-j)^2 + 4n(1-\alpha)\{n(1-\alpha)-(\vartheta-1)\}\frac{1}{((\vartheta-1)!)^2}\left(\prod_{j=0}^{\vartheta-2}(2n(1-\alpha)-j)\right)^2$$

$$= \frac{1}{((\vartheta-2)!)^2}\prod_{j=0}^{\vartheta-2}(2n(1-\alpha)-j)^2\left[1 + 4n(1-\alpha)\{n(1-\alpha)(-(\vartheta-1))\}\frac{1}{(\vartheta-1)^2}\right]$$

$$= \frac{1}{((\vartheta-1)!)^2}\prod_{j=0}^{\vartheta-2}(2n(1-\alpha)-j)^2\left[2n(1-\alpha)-(\vartheta-1)\right]^2 = \frac{1}{((\vartheta-1)!)^2}\prod_{j=0}^{\vartheta-1}(2n(1-\alpha)-j)^2.$$

The identity (2.3) therefore follows by induction on $\vartheta$. This completes the proof of Lemma 1.

Throughout in the sequel, let $I_k(n)$ denote the semiclosed interval $[\frac{k}{n}, \frac{k+1}{n})$, $k = 0,1,\ldots,n\text{-}1$.

**Lemma 2.** *Let the function* $f(z) = z + \sum_{n=2}^{\infty} a_n z^n$ *be in the class* $S^*(\alpha)$ *and* $a_{-n+\vartheta}^{(-n)}$ *be given by (2.2).*

*Then, for* $\alpha \in I_k(n)$,

$$\left|a_{-n+\vartheta}^{(-n)}\right| \leq \begin{cases} \dfrac{\Gamma(2n(1-\alpha)+1)}{\Gamma(\vartheta+1)\,\Gamma(2n(1-\alpha)+1-\vartheta)} & ,\vartheta = 1,\cdots\cdots,n-k \qquad (2.4) \\[2ex] \dfrac{\Gamma(2n(1-\alpha)+1)}{\vartheta\,\Gamma(n-k)\,\Gamma\big(2n(1-\alpha)+1-(n-k)\big)} & ,\vartheta = n-k+1,\cdots \qquad (2.5) \end{cases}$$

*In particular, if* $\alpha \in I_{n-1}(n)$, *then*

$$\left|a_{-n+\vartheta}^{(-n)}\right| \leq \frac{2n(1-\alpha)}{\vartheta}, \quad \vartheta = 1,2,\cdots \qquad (2.6)$$

*The estimates (2.4) and (2.6) are sharp.*

**Remark.** *By allowing k to vary from 0 to n-1, the estimates (2.4), (2.6) and (2.6) give estimates on the all the coefficients* $\left|a_{-n+\vartheta}^{(-n)}\right|$ *for a function* $f \in S^*(\alpha)$, $0 \leq \alpha < 1$.

**Proof.** By a direct calculation, $\big(z(1/(f(z))^n)'\big)/\big(-n(1/(f(z))^n)\big) = zf'(z)/f(z)$. Thus, $\big(z(1/(f(z))^n)'\big)/\big(-n(1/f(z)^n)\big) = (1+(1-2\alpha)w(z))/(1-w(z))$, $z \in D$, for a function $w(z) = \sum_{m=1}^{\infty} w_m z^m$ analytic in $D$ and satisfying the conditions of Schwarz Lemma. Equivalently,

$$z\big(1/(f(z))^n\big)' + \big(n/(f(z))^n\big) = \left[z\big(1/(f(z))^n\big)' - \big(n(1-2\alpha)/(f(z))^n\big)\right] w(z).$$ The Substitution of the corresponding series expansions of the functions in this identity and a simplification gives

$$\sum_{m=1}^{\infty} m\, a_{-n+m}^{(-n)} z^m = \left[-2n(1-\alpha) + \sum_{m=1}^{\infty} \{m - 2n(1-\alpha)\}\, a_{-n+m}^{(-n)} z^m\right] w(z) \qquad (2.7)$$

Equating coefficients on both sides of (2.7), it is observed that, for every $\vartheta = 1,2,\cdots$, the coefficient $a_{-n+\vartheta}^{(-n)}$ depends only on $a_{-n+1}^{(-n)}, a_{-n+2}^{(-n)}, \cdots, a_{-n+\vartheta-1}^{(-n)}$. Hence (2.7) can be rearranged as

$$\sum_{m=1}^{\vartheta} m\, a_{-n+m}^{(-n)} z^m + \sum_{m=\vartheta+1}^{\infty} b_m z^m = \left[-2n(1-\alpha) + \sum_{m=1}^{\vartheta-1} \{m - 2n(1-\alpha)\}\, a_{-n+m}^{(-n)} z^m\right] \sum_{p=1}^{\infty} w_p z^p, \quad \vartheta = 1,2,3,\cdots$$

the second sum in the left hand side being convergent in $D$. The inequality $|w(z)| < 1$ and Parseval's theorem give

$$\sum_{m=1}^{\vartheta} m^2 \left|a_{-n+m}^{(-n)}\right|^2 \leq 4n^2(1-\alpha)^2 + \sum_{m=1}^{\vartheta-1} \{m - 2n(1-\alpha)\}^2 \left|a_{-n+m}^{(-n)}\right|^2.$$

Equivalently,

$$\vartheta^2 \left|a_{-n+\vartheta}^{(-n)}\right|^2 \leq 4n(1-\alpha)\left[n(1-\alpha) + \sum_{m=1}^{\vartheta-1}\{n(1-\alpha)-m\}\left|a_{-n+m}^{(-n)}\right|^2\right] \quad . \tag{2.8}$$

The sign of each term inside the summation symbol on the right hand side of (2.8) depends on the sign of the expression $n(1-\alpha)-m, m=1,\cdots,\vartheta-1$. To determine the sign of this expression, we need to partition the interval $0 \leq \alpha < 1$ into $n$ semi-closed intervals $I_k(n), k = 0,1,\cdots,n-1$. For any fixed $k$, if $\alpha \in I_k(n), k=0,1,\cdots,n-1$, then $n-k-1 < n(1-\alpha) \leq n-k$ so that $n(1-\alpha)-m > 0$ if $m=1,\ldots,n-k-1$ and $n(1-\alpha)-m \leq 0$ if $m=n-k,\ldots$. Considering only nonnegative contributions in the right hand summation in (2.8), it follows by using the above inequalities that, for $\vartheta = 1,\ldots,n-k$,

$$\vartheta^2 \left|a_{-n+\vartheta}^{(-n)}\right|^2 \leq 4n(1-\alpha)\left[n(1-\alpha) + \sum_{m=1}^{\vartheta-1}\{n(1-\alpha)-m\}\left|a_{-n+m}^{(-n)}\right|^2\right] \tag{2.9}$$

while, if $\vartheta = n-k+1,\ldots,$ then

$$\vartheta^2 \left|a_{-n+\vartheta}^{(-n)}\right|^2 \leq 4n(1-\alpha)\left[n(1-\alpha) + \sum_{m=1}^{n-k-1}\{n(1-\alpha)-m\}\left|a_{-n+m}^{(-n)}\right|^2 + \sum_{m=n-k}^{\vartheta-1}\{n(1-\alpha)-m\}\left|a_{-n+m}^{(-n)}\right|^2\right.$$

$$\leq 4n(1-\alpha)\left[n(1-\alpha) + \sum_{m=1}^{n-k-1}\{n(1-\alpha)-m\}\left|a_{-n+m}^{(-n)}\right|^2\right] \tag{2.10}$$

We now use induction on $\vartheta$. For $\vartheta = 1$, it follows from (2.9) that $\left|a_{-n+1}^{(-n)}\right| \leq 2n(1-\alpha)$, giving the estimate (2.4) in this case. Now let, for $\vartheta = 1,2,\cdots,n-k-1$, the estimate

$$\left|a_{-n+\vartheta}^{(-n)}\right| \leq \prod_{j=0}^{\vartheta-1} \frac{2n(1-\alpha)-j}{j+1} \tag{2.11}$$

hold. Then, using (2.9), (2.11) and Lemma 1, it follows that, for $\vartheta = 1,\cdots,n-k$,

$$\vartheta^2 \left|a_{-n+\vartheta}^{(-n)}\right|^2 \leq 4n(1-\alpha)\left[n(1-\alpha) + \sum_{m=1}^{\vartheta-1}\{n(1-\alpha)-m\}\left(\prod_{j=0}^{m-1}\frac{2n(1-\alpha)-j}{j+1}\right)^2\right]$$

$$= \frac{1}{((\vartheta-1)!)^2}\prod_{j=0}^{\vartheta-1}(2n(1-\alpha)-j)^2 \tag{2.12}$$

Thus, for $\vartheta = 1,\cdots,n-k$,

$$\left|a_{-n+\vartheta}^{(-n)}\right| \leq \prod_{j=0}^{\vartheta-1} \frac{2n\ (1-\alpha)-j}{j+1} = \frac{\Gamma\left(2n\ (1-\alpha)+1\right)}{\Gamma(\vartheta+1)\ \Gamma\left(2n(1-\alpha)+1-\vartheta\right)} \qquad (2.13)$$

This establishes the inequality (2.4).

Next, if $\vartheta = n-k+1, n-k+2, \cdots\cdots$, using (2.10), the induction hypothesis (2.11) and Lemma 1, we get

$$\vartheta^2 \left|a_{-n+\vartheta}^{(-n)}\right|^2 \leq 4n\ (1-\alpha)\left[n\ (1-\alpha)+\sum_{m=1}^{n-k-1}\{n\ (1-\alpha)-m\}\left(\prod_{j=0}^{m-1}\frac{2n\ (1-\alpha)-j}{j+1}\right)^2\right]$$

$$= \frac{1}{\left((n-k-1)!\right)^2}\prod_{j=0}^{n-(k+1)}\left(2n(1-\alpha)-j\right)^2 = \left(\frac{\Gamma\left(2n\ (1-\alpha)+1\right)}{\Gamma(n-k)\ \Gamma\left(2n\ (1-\alpha)+1-(n-k)\right)}\right)^2$$

The above inequality yields the estimate (2.5).

For $k = n-1$, the estimates (2.4) and (2.5) respectively reduce to $\left|a_{-n+1}^{(-n)}\right| \leq 2n\ (1-\alpha)$ and $\left|a_{-n+\vartheta}^{(-n)}\right| \leq \frac{2n\ (1-\alpha)}{\vartheta}$ for $\vartheta = 2, 3, \cdots\cdots$. Combining the above inequalities, (2.6) follows.

Equality holds in (2.4) for every $\vartheta = 1, \cdots\cdots, n-k$ for $(-n)^{th}$ power of the function $K_\alpha(z)$ defined in (1.10). On the other hand for each $\vartheta = n-k+1, \ldots$, $(-n)^{th}$ power of the function $K_{\alpha,\vartheta}(z)$, defined in (1.10), provides the sharpness for the estimate (2.6). This completes the proof of Lemma 2.

We also need in the sequel the following result of Jabotinsky [8]:

**Lemma 3.** *If the function $f$, given by (2.1), is in $\Omega$ then $f^{-1} \in \Omega$. Further, if*

$$f^{-1}(w) = \sum_{n=1}^{\infty} A_n w^n$$

*then,*

$$A_n^{(p)} = \frac{p}{n}a_{-p}^{(-n)}, \quad n = 1, 2, \cdots; \quad p = \pm 1, \pm 2, \cdots \qquad (2.14)$$

*and $A_0^{(p)}$ is defined by*

$$\sum_{p=-\infty}^{\infty} A_0^{(p)} z^{-p-1} = \frac{f'(z)}{f(z)} \qquad . \qquad (2.15)$$

## 3. MAIN RESULTS

The following theorem gives the coefficient estimates for the inverse of a function in the class $S^*(\alpha)$:

**Theorem 1.** *Let* $f \in S^*(\alpha)$, $0 \leq \alpha < 1$ *and, for* $|w| < \frac{1}{4}$,

$$f^{-1}(w) = w + \sum_{n=2}^{\infty} A_n w^n \qquad (3.1)$$

*Then,*

(a) *for* $\alpha \in I_0(n) \cup I_1(n)$,

$$|A_n| \leq \frac{\Gamma(2n(1-\alpha)+1)}{\Gamma(n+1)\,\Gamma(2n(1-\alpha)+2-n)} \qquad (3.2)$$

(b) *for* $\alpha \in I_k(n)$, $k = 2, \cdots, n-2$,

$$|A_n| \leq \frac{\Gamma(2n(1-\alpha)+1)}{n(n-1)\,\Gamma(n-k)\,\Gamma(2n(1-\alpha)+1+k-n)} \qquad (3.3)$$

(c) *for* $\alpha \in I_{n-1}$,

$$|A_n| \leq \frac{2(1-\alpha)}{n-1} \qquad (3.4)$$

*where,* $I_k(n) \equiv [\frac{k}{n}, \frac{k+1}{n})$, $k = 0, 1, \ldots, n-1$. *The estimates (3.2) and (3.4) are sharp.*

**Proof.** It is known (see e.g. [12]) that $A_n = (1/2\pi i n) \int_{|z|=r} (1/(f(z))^n)\, dz = (1/n)\, a_{-1}^{(-n)}$, $0 < r < 1$. Therefore, it is sufficient to find suitable estimates for $\left|a_{-1}^{(-n)}\right|$. To this end, taking $\vartheta = n-1$ in Lemma 2, using the appropriate inequality (2.4), (2.5) or (2.6) for different values of $k = 0, 1, \ldots, n-1$ and observing that only for $k=0$ and $k=1$ the inequality (2.4) is applicable, the following estimate is obtained for $\alpha \in [0, \frac{2}{n})$,

$$\left|a_{-1}^{(-n)}\right| \leq \frac{\Gamma(2n(1-\alpha)+1)}{\Gamma(n)\,\Gamma(2n(1-\alpha)+2-n)} \qquad (3.5)$$

which gives (a). Similarly, for $\alpha \in I_k(n)$, $k = 2, \ldots, n-2$, the inequality (2.5) yields

$$\left|a_{-1}^{(-n)}\right| \leq \frac{\Gamma(2n(1-\alpha)+1)}{(n-1)\,\Gamma(n-k)\,\Gamma(2n(1-\alpha)+1+k-n)}. \qquad (3.6)$$

This gives (b). Finally, for $\alpha \in I_{n-1}(n)$, the inequality (2.6) gives

$$\left|a_{-1}^{(-n)}\right| \leq \frac{2n\,(1-\alpha)}{n-1}. \tag{3.7}$$

Consequently, (c) follows.

It is easily verified that equality holds in (3.5) and (3.7) for the $(-n)^{th}$ power of the function $K_n(z)$ and $(-n)^{th}$ power of the function $K_{\alpha,n-1}(z)$ respectively. Thus, the estimates (3.2) and (3.4) are sharp. This completes the proof of Theorem 1.

**Remark.** *The sharp coefficient bounds of Krzysz, Libera and Zlotkiewicz [12] for $|A_2|$ and $|A_3|$ follow as a particular case of Theorem 1.*

**Remark:** *The sharp coefficient bounds of Krzyz, Libera and Zlotkiewicz [12] for $|A_2|$ and $|A_3|$ follow as a particular case of Theorem 1.*

The sharp coefficient estimates for functions in $\Sigma^*(\alpha)$, $0 \leq \alpha < 1$, are described by the following:

**Theorem 2.** *Let the function $g \in \Sigma^*(\alpha)$, $0 \leq \alpha < 1$, be given by the series*

$$g(z) = z + b_0 + \frac{b_1}{z} + \frac{b_2}{z^2} + \cdots\cdots, \quad z \in V.$$

*Then,*

$$|b_m| \leq \frac{2\,(1-\alpha)}{m+1}, \quad m = 0,1\ldots \tag{3.8}$$

*The estimate (3.8) is sharp.*

**Proof.** The mapping $f(z) \to g(z) = 1/f(1/z)$ establishes a one-to-one correspondence between $S^*(\alpha)$ and $\Sigma^*(\alpha)$. Since, $(zg'(z)/g(z)) = \left(z(1/f(1/z))'\right)/(1/f(1/z)) = ((1/z)f'(1/z))/f(1/z)$, this mapping too is one-to-one from $S^*(\alpha)$ onto $\Sigma^*(\alpha)$. We note that the coefficient expansion of $1/f(z)$ around origin is same as the coefficient expansion of $g(z)$ about infinity. Therefore,

$$\max_{g \in \Sigma^*(\alpha)} |b_m| = \max_{f \in S^*(\alpha)} \left|a_m^{(-1)}\right|, \quad m = 0,1\ldots \tag{3.9}$$

Thus, by Lemma 2,

$$\left|a_{-1+\vartheta}^{(-1)}\right| \leq \frac{2\,(1-\alpha)}{\vartheta}, \quad 0 \leq \alpha < 1;\ \vartheta = 1,\ldots \tag{3.10}$$

The inequality (3.10) can be equivalently expressed as

$$\left|a_m^{(-1)}\right| \leq \frac{2(1-\alpha)}{m+1}, \quad 0 \leq \alpha < 1;\ m = 0,1,\ldots \tag{3.11}$$

and the inequality (3.11) together with (3.9) gives (3. 8). It is easily seen that the function $g(z) = z\left(1-(1/z^{m+1})\right)^{2(1-\alpha)/(m+1)}$ belongs to the class $\Sigma^*(\alpha)$ and its $m^{th}$ coefficient equals $2(1-\alpha)/(m+1)$. Therefore, the estimate (3.8) is sharp. This completes the proof of Theorem 2.

The coefficient estimates for the inverse of a function in the class $\Sigma^*(\alpha)$ are found in the following:

**Theorem 3.** *Let the function $g \in \Sigma^*(\alpha), 0 \leq \alpha < 1$, have the series expansion $g^{-1}(w) = w + \sum_{n=0}^{\infty} B_n w^{-n}$ in some neighbourhood of infinity. Then,*

(a) $\quad |B_0| \leq 2\,(1-\alpha), \quad 0 \leq \alpha < 1$ \hfill (3.12)

(b) *For $\alpha \in I_k(n), \quad k = 0,\cdots\cdots n-2$,*
$$|B_n| \leq \frac{\Gamma(2n\,(1-\alpha)+1)}{n(n+1)\,\Gamma(n-k)\,\Gamma(2n\,(1-\alpha)+1-(n-k))} \quad (3.13)$$

(c) *For $\alpha \in I_{n-1}(n)$,*
$$|B_n| \leq \frac{2\,(1-\alpha)}{n+1} \quad (3.14)$$

*where, $I_k(n) \equiv [\frac{k}{n}, \frac{k+1}{n})$, $k = 0,1,...,n-1$. The estimates (3.12) and (3.14) are sharp.*

**Proof.** For any $g \in \Sigma^*(\alpha)$, $0 \leq \alpha < 1$, there exists $f \in S^*(\alpha)$ such that $g(z) = 1/f(1/z)$. It can be easily verified that $g^{-1}(w) = 1/f^{-1}(1/w)$. Since the coefficients in the expansion of $1/f^{-1}(w)$ around origin and those of $1/f^{-1}(1/w)$ about infinity are the same, we have
$$B_n = A_n^{(-1)}, \quad n = 0,1,2,\cdots. \quad (3.15)$$
By (2.15), with $\operatorname{Re} Q(z) > \alpha$, $Q(0) = 1$ and $z \in D$,
$$\sum_{p=-\infty}^{\infty} A_0^{(p)} z^{-p-1} = \frac{f'(z)}{f(z)} = \frac{1}{z} Q(z) \equiv \frac{1}{z}[1 + \sum_{n=1}^{\infty} q_n z^n]$$
Now, using the well known sharp estimate
$$|q_n| \leq 2\,(1-\alpha), \quad n = 1, 2, \cdots \quad (3.16)$$
we get,
$$\max_{g^{-1} \in \Sigma^{*-1}(\alpha)} |B_0| = \max_{f^{-1} \in S^{*-1}(\alpha)} |A_0^{-1}| = \operatorname{Max}|q_1| \leq 2\,(1-\alpha),$$

This establishes (3.12). Since the bound in (3.16) is sharp, it follows that the estimate (3.12) is also sharp.

For $n=1, 2,...$, (3.15) together with (2.14) gives

$$\max_{g^{-1} \in \Sigma^{*-1}(\alpha)} |B_n| = \max_{f^{-1} \in S^{*-1}(\alpha)} |A_n^{(-1)}| = (1/n) \max_{f^{-1} \in S^*(\alpha)} |a_1^{(-n)}|$$
(3.17)

It is easily verified that $a_1^{(-n)} = a_{-n+(n+1)}^{(-n)}$. Since the sharp estimate (2.4) in Lemma 2 is not applicable for any value of $k = 0,\cdots,n-1$, in order to get (3.13), the estimate (2.5) with $\vartheta = n+1$ has to be used for $k = 0,1,...,n-2$. This gives

$$|a_1^{(-n)}| \leq \frac{\Gamma(2n(1-\alpha)+1)}{(n+1)\,\Gamma(n-k)\,\Gamma(2n(1-\alpha)+1-(n-k))}$$
(3.18)

The estimate (3.13) now easily follows from (3.17) and (3.18). Similarly, the estimate (2.6) with $\vartheta = n+1$ gives,

$$|a_1^{(-n)}| \leq \frac{2n(1-\alpha)}{n+1}$$
(3.19)

By combining (3.17) and (3.19), the estimate (3.14) follows. Since the inequality (2.6) is sharp, it follows that the estimate (3.14) is sharp. This completes the proof of Theorem 3.

**Remark.** *The construction of a suitable example to exhibit the sharpness of inequality (2.5) seems to be quite involved and the sharpness of estimates (3.3) and (3.13) depend on the sharpness of the inequality (2.5).*